\documentclass{amsart}
\usepackage{amsmath,amssymb}
\usepackage[dvips]{graphics}
\usepackage[all]{xy}
\newtheorem{Theorem}{Theorem}[section]
\newtheorem{Lemma}[Theorem]{Lemma}
\newtheorem{Proposition}[Theorem]{Proposition}

\newtheorem{Example}[Theorem]{Example}

\newtheorem{Remark}[Theorem]{Remark}

\makeatletter
\@addtoreset{figure}{section}
\def\@thmcountersep{-}
\makeatother

\numberwithin{equation}{section}

\begin{document}

\title{Site-specific Gordian distances of spatial graphs}

\author{Kouki Taniyama}
\address{Department of Mathematics, School of Education, Waseda University, Nishi-Waseda 1-6-1, Shinjuku-ku, Tokyo, 169-8050, Japan}
\email{taniyama@waseda.jp}
\thanks{The author was partially supported by Grant-in-Aid for Challenging Exploratory Research (No. 15K13439) and Grant-in-Aid for Scientific Research(A) (No. 16H02145) , Japan Society for the Promotion of Science.}

\subjclass[2010]{Primary 57M25; Secondly 57M15.}

\date{}

\dedicatory{}

\keywords{knot, link, spatial graph, puzzle ring, site-specific Gordian distance, Milnor link, covering space}

\begin{abstract}
A site-specific Gordian distance between two spatial embeddings of an abstract graph is the minimal number of crossing changes from one to another where each crossing change is performed between two previously specified abstract edges of the graph. It is infinite in some cases. We determine the site-specific Gordian distance between two spatial embeddings of an abstract graph in certain cases. It has an application to puzzle ring problem. The site-specific Gordian distances between Milnor links and trivial links are determined. We use covering space theory for the proofs. 
\end{abstract}

\maketitle

\section{Introduction} 

Throughout this paper we work in the piecewise linear category. 
Let $G$ be a finite graph. We denote the set of all vertices of $G$ by $V(G)$ and the set of all edges of $G$ by $E(G)$. 
Let ${\mathcal E}=(E(G)\times E(G))/\sim$ be the quotient set of $E(G)\times E(G)$ under the equivalence relation $\sim$ defined by $(x,y)\sim(x,y)$ and $(x,y)\sim(y,x)$ for every $(x,y)\in E(G)\times E(G)$. We denote the equivalence class of $(x,y)$ by $[x,y]$. 
Let $f:G\to {\mathbb S}^3$ be an embedding of $G$ into the $3$-sphere ${\mathbb S}^3$. 
The embedding $f$ is said to be a {\it spatial embedding} of $G$. 
Let $x$ and $y$ be two edges of $G$. 
By an {\it $x$-$y$ crossing change on $f$} we mean a crossing change between a part of  $f(x)$ and a part of $f(y)$. 
When $x=y$, an $x$-$y$ crossing change is said to be a {\it self-crossing change on $x$}. 
Let ${\mathcal F}$ be a subset of ${\mathcal E}$. 
We say that a crossing change ${\rm C}$ belongs to ${\mathcal F}$ if ${\rm C}$ is an $x$-$y$ crossing change and $[x,y]$ is an element of ${\mathcal F}$. 
Let $f:G\to {\mathbb S}^3$ and $g:G\to {\mathbb S}^3$ be two spatial embeddings of $G$. We consider these embeddings up to ambient isotopy in ${\mathbb S}^3$. 
The {\it ${\mathcal F}$-Gordian distance between $f$ and $g$} is defined to be the minimal number of crossing changes each of which belongs to ${\mathcal F}$ that deforms $f$ into $g$. We denote it by $d_{\mathcal F}(f,g)$. 
In case that $f$ and $g$ are not transformed into each other by such crossing changes, $d_{\mathcal F}(f,g)=\infty$. 

Let $G$ be a planar graph. An embedding $t:G\to {\mathbb S}^3$ is said to be {\it trivial} if $t(G)$ is contained in a $2$-sphere in ${\mathbb S}^3$. 
Note that trivial embeddings are unique up to ambient isotopy \cite{Mason}. 
The image of a trivial embedding is also said to be {\it trivial}. 
Let $f:G\to {\mathbb S}^3$ be a spatial embedding. The ${\mathcal F}$-unknotting number of $f$, denoted by $u_{\mathcal F}(f)$, is defined to be the ${\mathcal F}$-Gordian distance between $f$ and $t$. Namely $u_{\mathcal F}(f)=d_{\mathcal F}(f,t)$. 

In this paper we mainly concern with the case that ${\mathcal F}$ is a singleton. Namely ${\mathcal F}=\{[x,y]\}$ for some edges $x$ and $y$ of $G$. 
When $E(G)=\{e_1,\cdots,e_n\}$ we denote $e_i$-$e_j$ crossing change by $i$-$j$ crossing change and $d_{\{[e_i,e_j]\}}(f,g)$ by $d_{i,j}(f,g)$ for the simplicity. 
When $G$ is planar $u_{\{[e_i,e_j]\}}(f)$ is denoted by $u_{i,j}(f)$. 

Let $\mu$ be a natural number and $L_{\mu}$ a graph with $\mu$ vertices $v_1,\cdots,v_\mu$ and $\mu$ loops $\ell_1,\cdots,\ell_\mu$ such that $v_i$ and $\ell_i$ are incident for each $i\in\{1,\cdots,\mu\}$. Then a spatial embedding of $L_\mu$ can be regarded as an ordered oriented link in ${\mathbb S}^3$. In this paper we ignore a vertex of a graph incident to exactly one edge that is a loop so long as no confusion occurs. 

\vskip 3mm

\begin{figure}[htbp]
      \begin{center}
\scalebox{0.5}{\includegraphics*{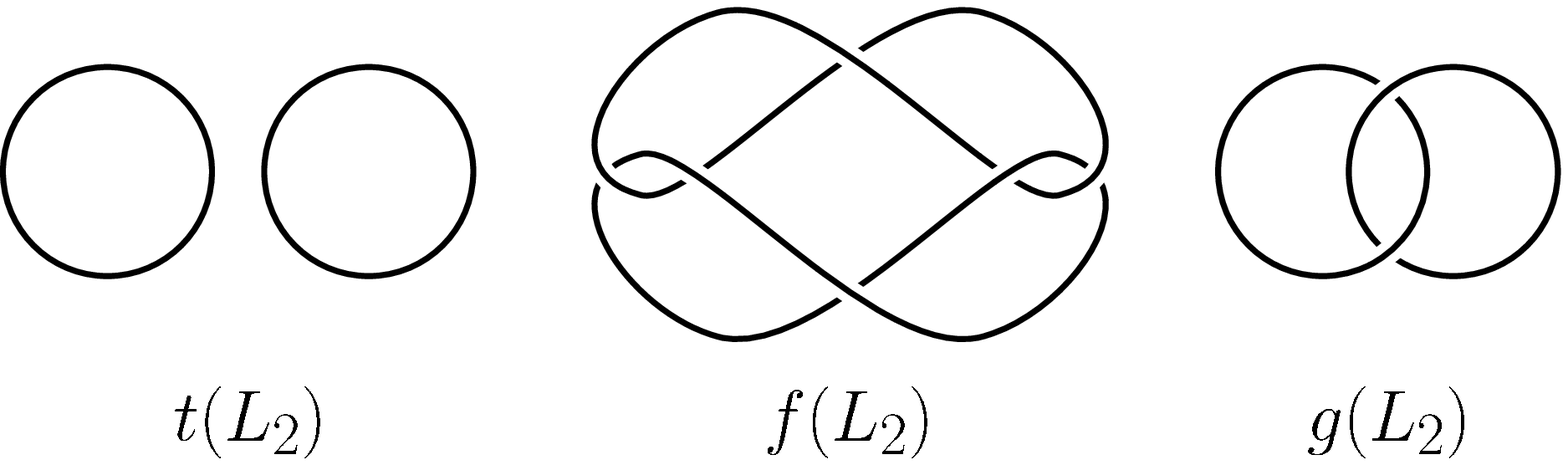}}
      \end{center}
   \caption{}
  \label{two-component-links}
\end{figure}

\begin{Example}\label{two-component link}
{\rm
Let $t:L_2\to {\mathbb S}^3$, $f:L_2\to {\mathbb S}^3$ and $g:L_2\to {\mathbb S}^3$ be spatial embeddings as illustrated in Figure \ref{two-component-links}. The orders and orientations are arbitrary. Then we have $u_{1,1}(f)=d_{1,1}(f,t)=u_{2,2}(f)=d_{2,2}(f,t)=1$, $u_{1,2}(f)=d_{1,2}(f,t)=2$, $u_{1,1}(g)=d_{1,1}(g,t)=u_{2,2}(g)=d_{2,2}(g,t)=\infty$, $u_{1,2}(g)=d_{1,2}(g,t)=1$, $d_{1,1}(f,g)=d_{2,2}(f,g)=\infty$ and $d_{1,2}(f,g)=1$. 
}
\end{Example}

%

\begin{figure}[htbp]
      \begin{center}
\scalebox{0.65}{\includegraphics*{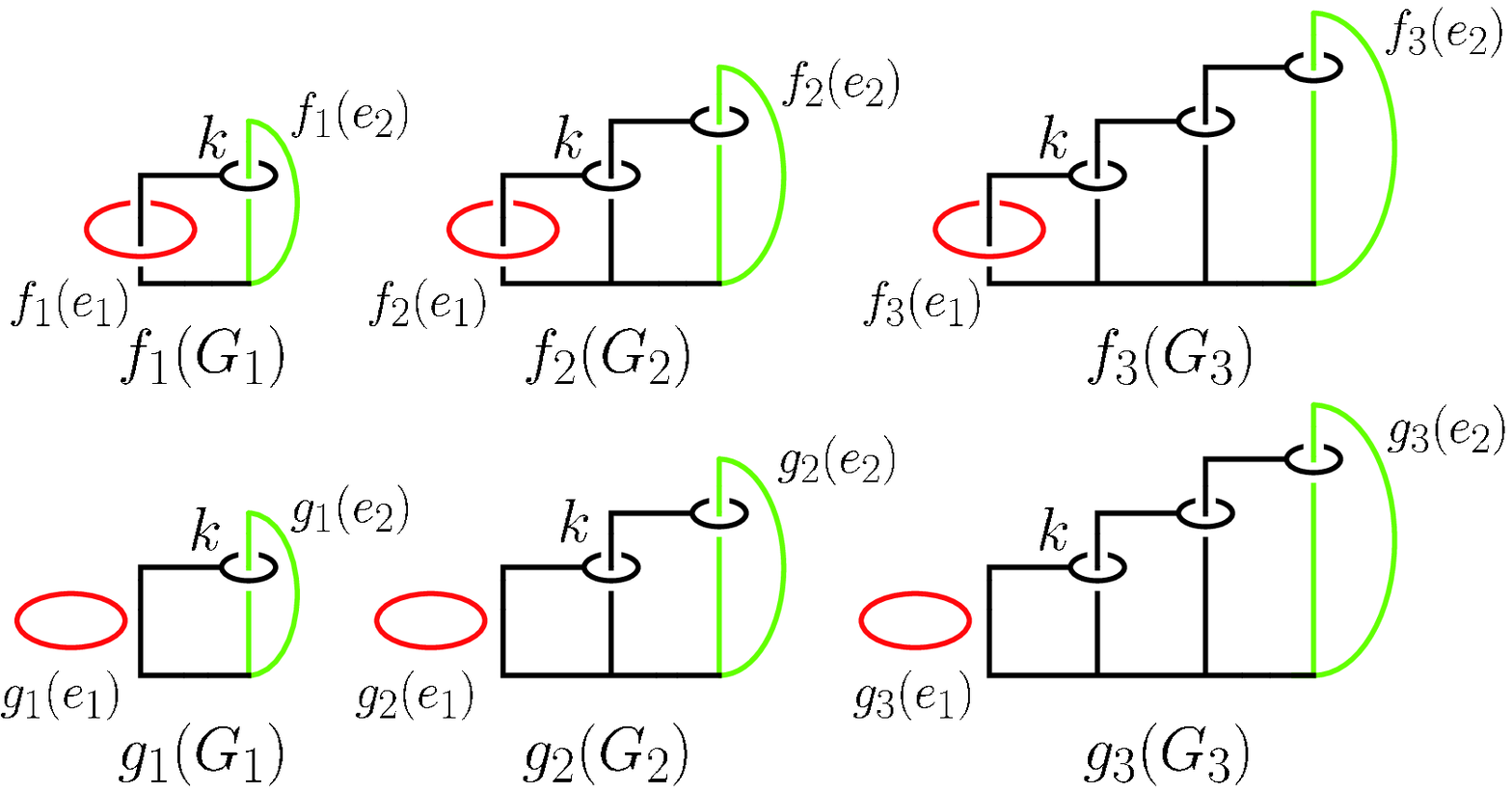}}
      \end{center}
   \caption{}
  \label{puzzle1}
\end{figure}

Let $n$ be a natural number. Let $G_n$ be an abstract graph, $e_1$ and $e_2$ abstract loops of $G_n$, and $f_n:G_n\to {\mathbb S}^3$ and $g_n:G_n\to {\mathbb S}^3$ spatial embeddings of $G_n$ as illustrated in Figure \ref{puzzle1}. Here $f_n$ and $g_n$ differ only on $e_1$. 
In Figure \ref{puzzle1} only the cases $n=1,2,3$ are illustrated. The general case $n$ is an obvious generalization of these cases. 

Then the following theorem is shown in \cite{P-S}. 

\vskip 3mm

\begin{Theorem}\label{puzzle-1}\cite{P-S}
Let $n$ be a natural number. Let $f_n:G_n\to {\mathbb S}^3$ and $g_n:G_n\to {\mathbb S}^3$ be spatial embeddings of $G_n$ as illustrated in Figure \ref{puzzle1}. 
Then $d_{1,2}(f_n,g_n)=2^n$. 
\end{Theorem}

\vskip 3mm

\begin{figure}[htbp]
      \begin{center}
\scalebox{0.6}{\includegraphics*{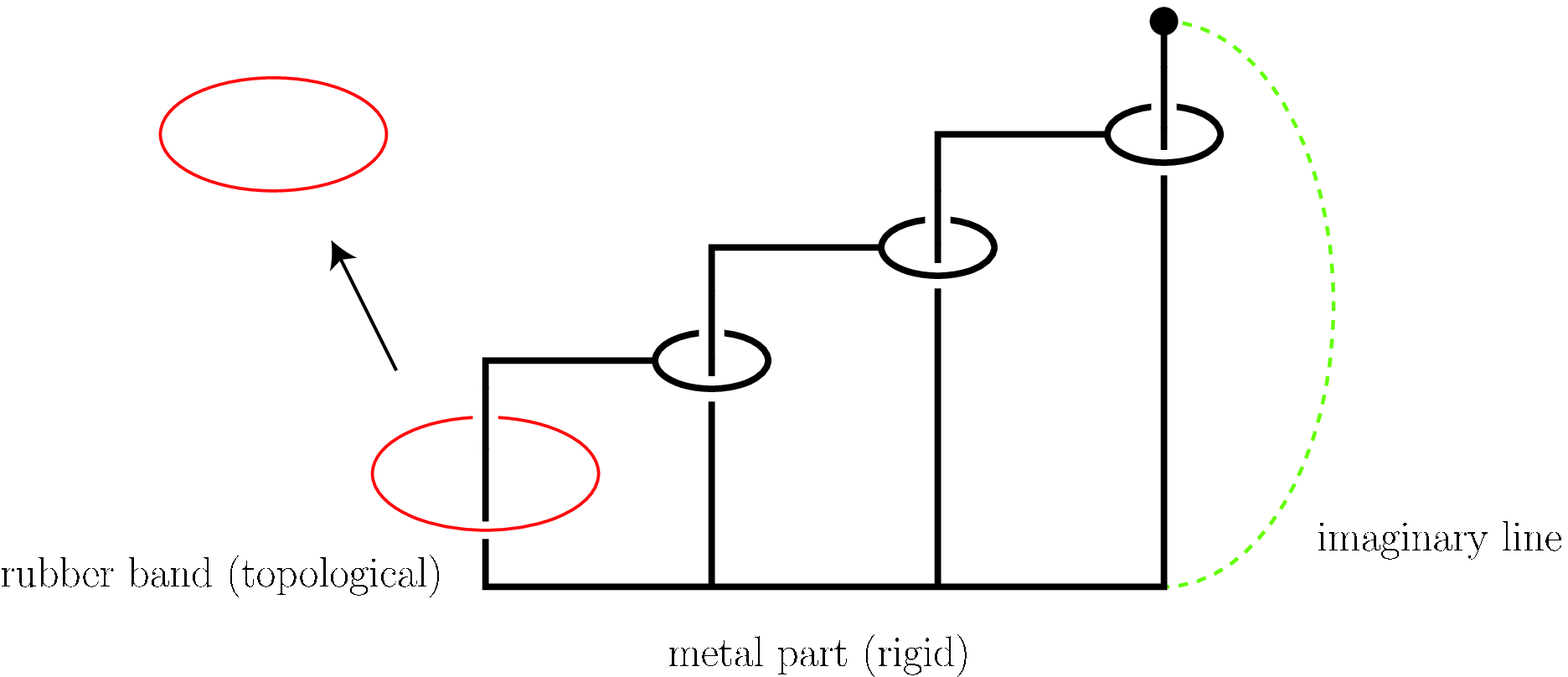}}
      \end{center}
   \caption{}
  \label{puzzle3}
\end{figure}

Note that Theorem \ref{puzzle-1} is originally an answer to a question on puzzle ring raised by Kauffman in \cite{Kauffman}. 
Suppose that $f_n(e_1)$ is a rubber band, the right half of $f_n(e_2)$ is not real but imaginary, and the rest of the image of $f_n(G_n)$ is made of metal. 
Here we think that the rubber band is topological and the metal part is rigid. Then it follows by Theorem \ref{puzzle-1} that the rubber band must go across the imaginary part $2^n$ times before it comes away from the metal part. See Figure \ref{puzzle3} where the case $n=3$ is illustrated. See \cite{P-S} for more detail. 

The proof of Theorem \ref{puzzle-1} in \cite{P-S} is based on a group-theoretical argument. In this paper we give a new proof of Theorem \ref{puzzle-1} based on covering space theory. The method of proof is applicable to other similar problems. As examples we  show Theorem \ref{puzzle-2} and Theorem \ref{Milnor} below. So far as the author's understanding the argument in \cite{P-S} does not work for the proof of Theorem \ref{puzzle-2}. 

\begin{figure}[htbp]
      \begin{center}
\scalebox{0.6}{\includegraphics*{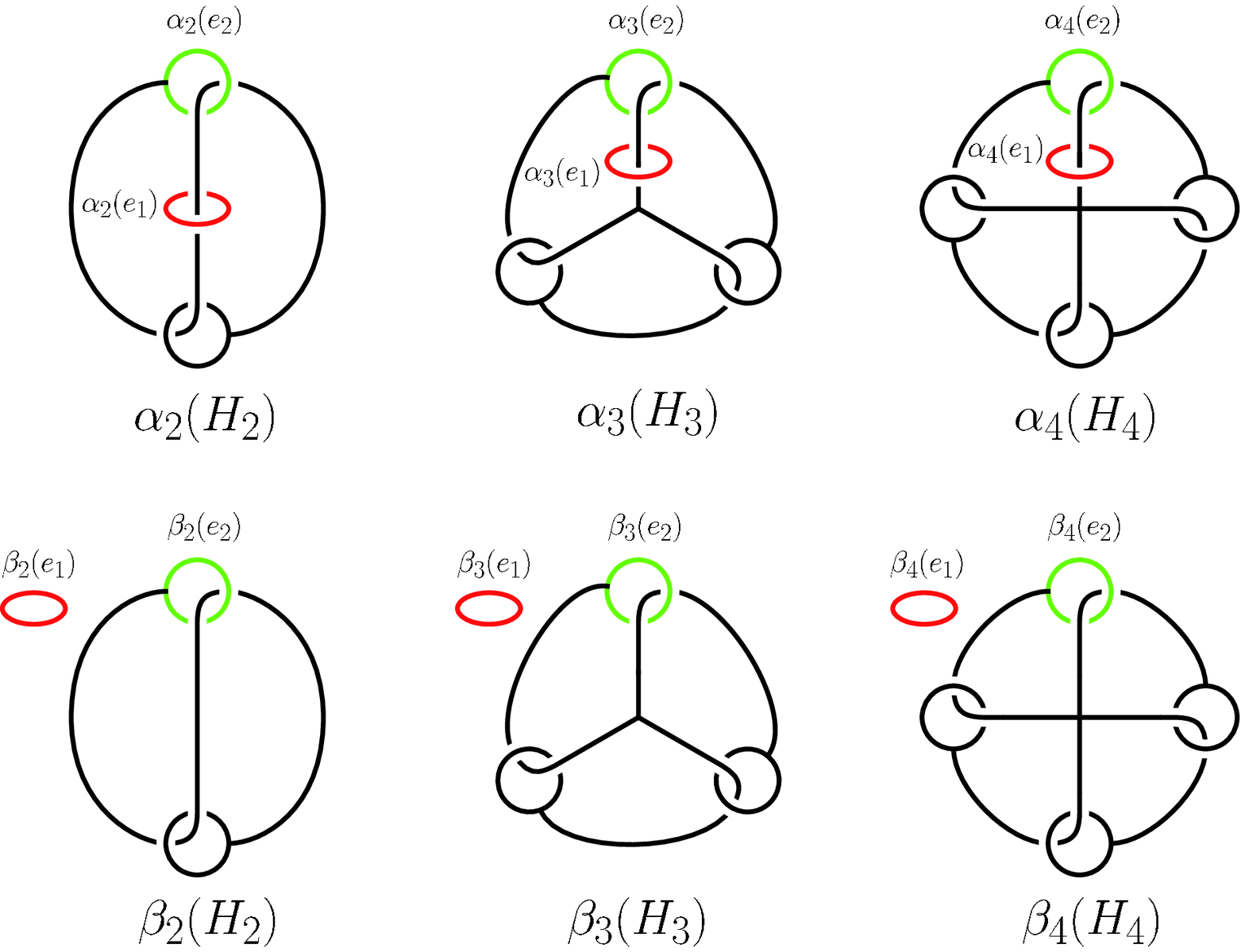}}
      \end{center}
   \caption{}
  \label{puzzle2}
\end{figure}

Let $n$ be a natural number with $n\geq2$. Let $H_n$ be an abstract graph, $e_1$ and $e_2$ abstract loops of $H_n$, and $\alpha_n:H_n\to {\mathbb S}^3$ and $\beta_n:H_n\to {\mathbb S}^3$ spatial embeddings of $H_n$ as illustrated in Figure \ref{puzzle2}. Here $\alpha_n$ and $\beta_n$ differ only on $e_1$. 
In Figure \ref{puzzle2} only the cases $n=2,3,4$ are illustrated. The general case $n$ is an obvious generalization of these cases. 
Then we have the following theorem. 

\begin{Theorem}\label{puzzle-2}
Let $n$ be a natural number with $n\geq2$. Let $\alpha_n:H_n\to {\mathbb S}^3$ and $\beta_n:H_n\to {\mathbb S}^3$ be spatial embeddings of $H_n$ as illustrated in Figure \ref{puzzle2}. 
Then $d_{1,2}(\alpha_n,\beta_n)=2^n-2$. 
\end{Theorem}

\vskip 3mm

Let $J_n$ be the connected component of $H_n$ containing $e_2$. 
Then $\alpha_n(J_n)$ is nontrivial \cite{Taniyama2}. Since every proper subgraph of $\alpha_n(J_n)$ is trivial we see that the fundamental group of the complement ${\mathbb S}^3\setminus\alpha_n(J_n)$ is not a free group \cite{ST}. Therefore the argument in \cite{P-S} regarding $\alpha_n(e_1)$ as an element of the fundamental group will not work here. 
On the other hand it is not hard to see, for example by repeated use of the replacement of crossing changes illustrated in Figure \ref{replacing} and Figure \ref{replacing2}, that $\alpha_n(H_n)$ is deformed into a trivial spatial graph by self-crossing changes on the edge $e_2$. 
Therefore the argument in \cite{P-S} regarding $\alpha_n(e_2)$ as an element of the fundamental group will not work, too. 
Therefore Theorem \ref{puzzle-2} is an entirely new result to the best of the author's understanding. 

\begin{figure}[htbp]
      \begin{center}
\scalebox{0.5}{\includegraphics*{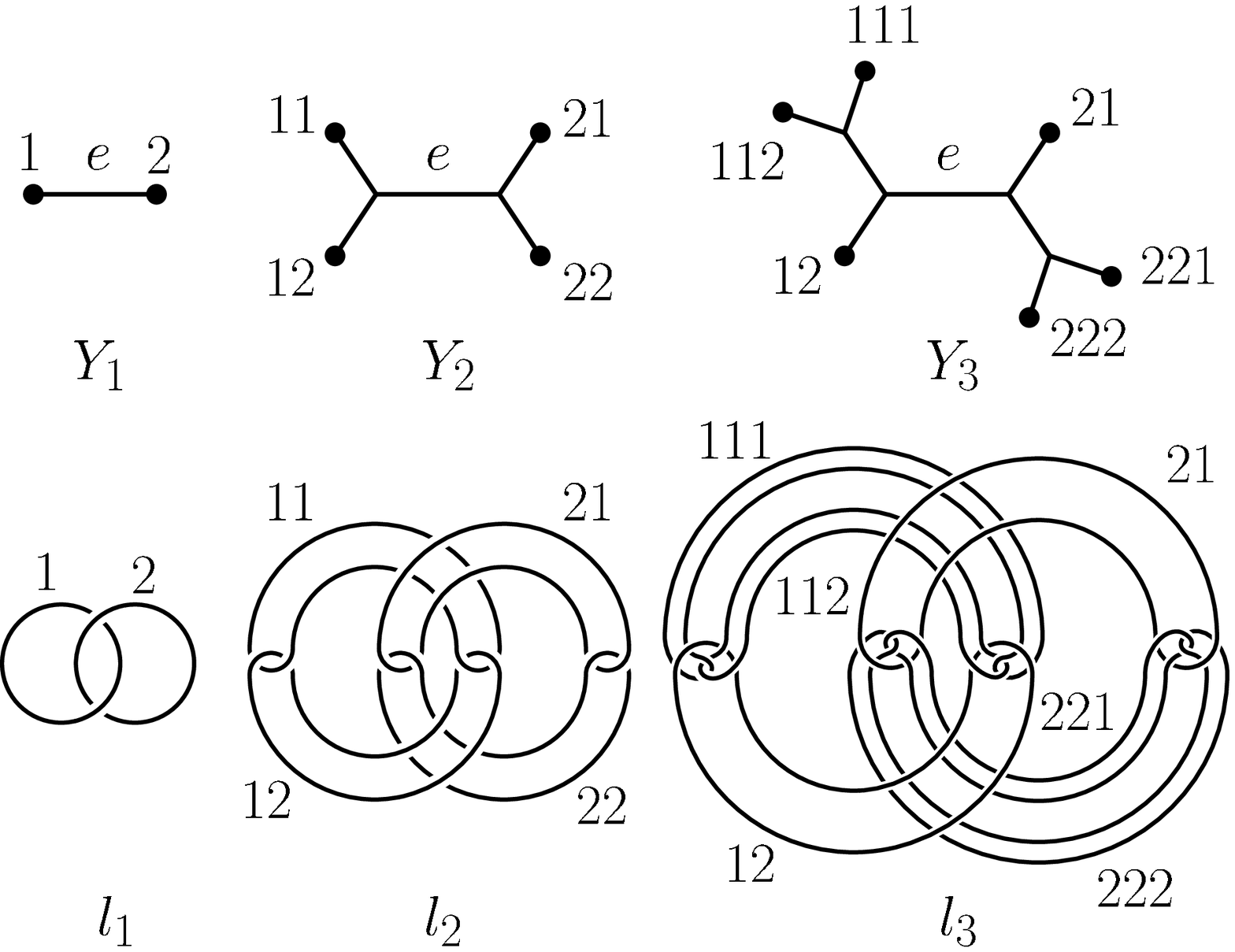}}
      \end{center}
   \caption{}
  \label{Milnor-type-link}
\end{figure}

Let $T$ be a uni-trivalent tree. Namely $T$ is a finite connected graph without cycles and $V(T)=V_1(T)\cup V_3(T)$ where $V_i(G)$ denotes the set of all degree $i$ vertices of a graph $G$. Let $e$ be an edge of $T$. Let $Y_1$ be the minimal subgraph of $T$ containing the edge $e$. Let $Y_2$ be the minimal subgraph of $T$ containing $Y_1$ and all edges of $T$ adjacent to $e$. Suppose that $Y_k$ is defined. Then $Y_{k+1}$ is defined to be the minimal subgraph of $T$ containing $Y_k$ and all edges of $T$ adjacent to some edge of $Y_k$. Then there are a natural number $n$ and a sequence of subgraphs $Y_1, Y_2, \cdots, Y_n$ of $T$ such that $Y_1\subset Y_2\subset\cdots\subset Y_n=T$ and $Y_{n-1}$ is a proper subgraph of $T$ if $n\geq2$. For each $Y_i$ we assign an unoriented link $l_i\subset{\mathbb S}^3$ and a bijection from $V_1(Y_i)$ to the set of the components of $l_i$, denoted by $C(l_i)$, as follows. 
For $Y_1$ we assign a Hopf link $l_1$. We assign an arbitrary bijection from $V_1(Y_1)$ to $C(l_1)$. Now suppose that for $Y_k$ a link $l_k$ and a bijection from $V_1(Y_k)$ to $C(l_k)$ are assigned. For each vertex $v\in V_1(Y_k)\cap V_3(Y_{k+1})$ we replace the component of $l_k$ corresponding to $v$ by a pair of components forming a Bing double \cite{Bing}. The two vertices in $V(Y_{k+1})\setminus V(Y_k)$ adjacent  to $v$ correspond to the pair of components forming the Bing double. Other correspondences are inherited from that of $Y_k$ and $l_k$. See for example Figure \ref{Milnor-type-link} where the sequences of numbers describe the correspondence between $V_1(Y_k)$ and the components of $l_k$. We call $L_T=l_n$ the unoriented link associated with $T$ and $e$ the defining edge of $L_T$. 
Let $v$ and $w$ be two vertices of $T$. Let $P(v,w)$ be the path in $T$ with $\partial P(v,w)=\{v,w\}$. Let $d(v,w)$ be the number of the edges of $P(v,w)$. Let $a$ and $b$ be mutually distinct components of $L_T$. Let $v$ and $w$ be the corresponding vertices in $V_1(T)$. Set $d(a,b)=d(v,w)$.

\begin{Remark}\label{clasper}
{\rm
(1) We will see in Proposition \ref{well-definedness} that the link type of $l_n$ is independent of the choice of the defining edge $e$ of $T$. We also note that $l_n$ is obtained from a trivial link by applying a $C_m$-move along a simple tree clasper whose underlying tree is $T$ \cite{Habiro}. In this sense the well-definedness is already known. 

\noindent
(2) Let $T$ be a uni-trivalent tree such that every element of $V_3(T)$ is adjacent to an element of $V_1(T)$. Then $L_T$ is a Milnor link first appeared in \cite{Milnor}. 

\noindent
(3) For an unoriented link $L$ and a subset ${\mathcal F}$ of ${\mathcal E}$, the ${\mathcal F}$-unknotting number $u_{\mathcal F}(L)$ of $L$ is well-defined since any two oriented trivial links of the same number of components are mutually ambient isotopic.  
}
\end{Remark}

\vskip 3mm

\begin{Theorem}\label{Milnor}
Let $T$ be a uni-trivalent tree. Let $L_T$ be the unoriented link associated with $T$. Let $a$ and $b$ be mutually distinct components of $L_T$ and ${\mathcal F}=\{[a,b]\}$. Then $u_{\mathcal F}(L_T)=2^{d(a,b)-1}$. 
\end{Theorem}

\vskip 3mm

\section{Proofs}\label{section2} 

\begin{figure}[htbp]
      \begin{center}
\scalebox{0.55}{\includegraphics*{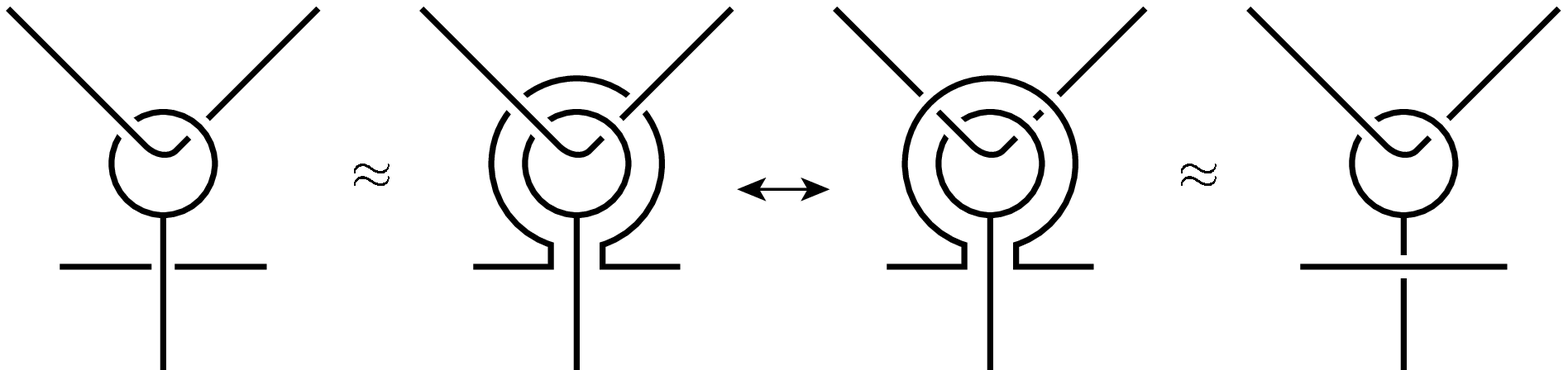}}
      \end{center}
   \caption{}
  \label{replacing}
\end{figure}

\begin{figure}[htbp]
      \begin{center}
\scalebox{0.65}{\includegraphics*{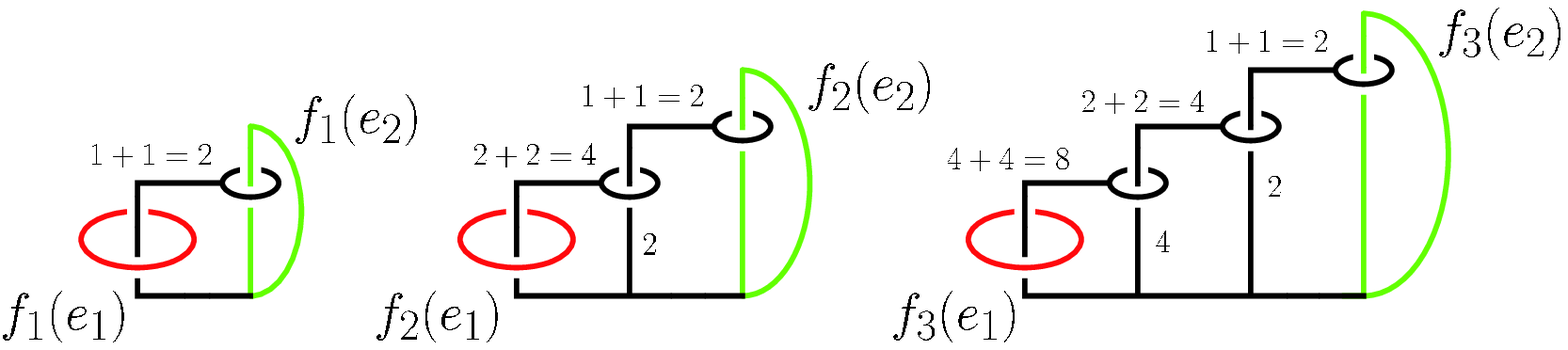}}
      \end{center}
   \caption{}
  \label{puzzle1-1}
\end{figure}

\noindent{\bf Proof of Theorem \ref{puzzle1}.} 
First we give a proof of $d_{1,2}(f_n,g_n)\leq 2^n$. The proof is essentially the same as that in \cite{P-S} or \cite{Taniyama}. 
We repeatedly use the technic of replacing crossing changes as illustrated in Figure \ref{replacing}. 
In Figure \ref{puzzle1-1} the number attached to an edge describes the number of crossing changes with $f_n(e_2)$ that realize a crossing change with the edge. Then we see $d_{1,2}(f_n,g_n)\leq 2^n$. 

\begin{figure}[htbp]
      \begin{center}
\scalebox{0.55}{\includegraphics*{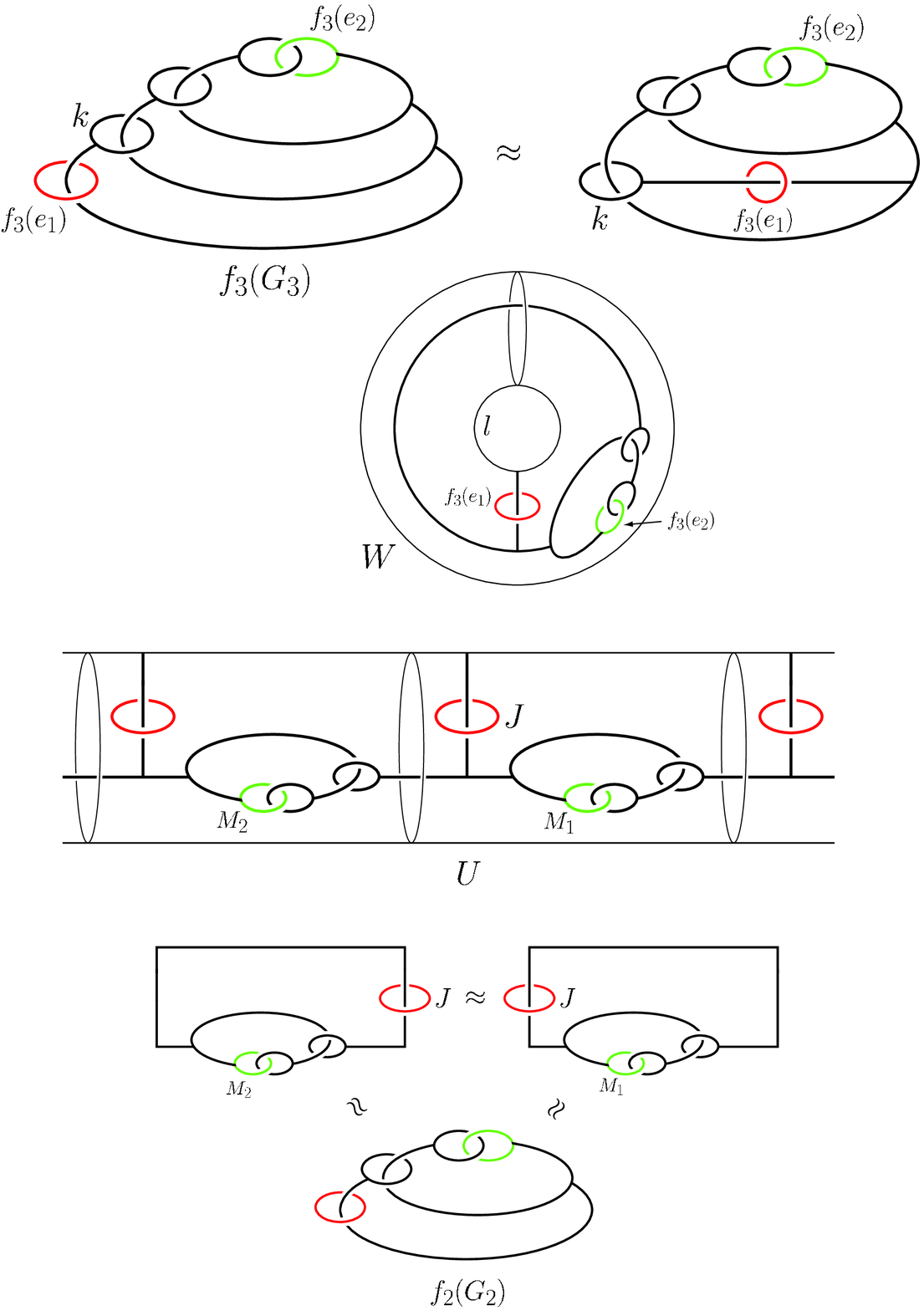}}
      \end{center}
   \caption{}
  \label{covering1}
\end{figure}

\begin{figure}[htbp]
      \begin{center}
\scalebox{0.55}{\includegraphics*{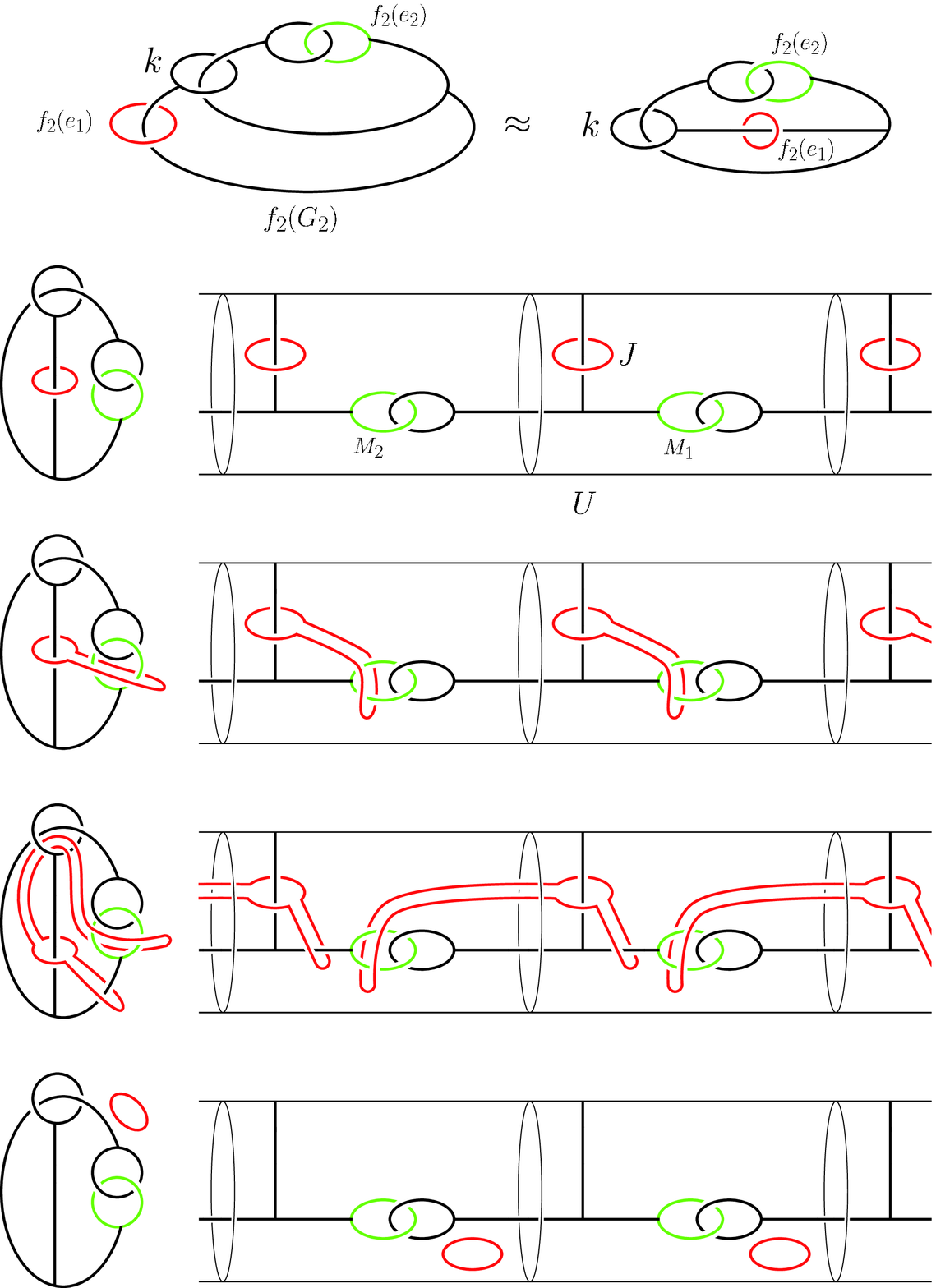}}
      \end{center}
   \caption{}
  \label{covering2}
\end{figure}

Next we give a new proof of $d_{1,2}(f_n,g_n)\geq 2^n$ by an induction on $n$. 
Let $k\subset f_n(G_n)\cap g_n(G_n)$ be the trivial knot illustrated in Figure \ref{puzzle1}. Let $N(k)$ be a regular neighbourhood of $k$ in ${\mathbb S}^3$. We may suppose that $f_n(G_n)\cap N(k)=g_n(G_n)\cap N(k)$ is a simple arc. Since $k$ is trivial $W={\mathbb S}^3\setminus{\rm int}N(k)$ is a solid torus. Let $\varphi:U\to W$ be the universal covering projection. Then we have the preimage $\varphi^{-1}(f_n(G_n)\cap W)$ in $U$. 
See Figure \ref{covering1} where the case $n=3$ is illustrated. 
Since $U$ is homeomorphic to ${\mathbb D}^2\times{\mathbb R}$, ${\mathbb D}^2\times{\mathbb R}$ is a subset of ${\mathbb R}^3$ and ${\mathbb R}^3$ is a subset of its one-point compactification ${\mathbb S}^3$, we may suppose that $U$ is a subset of ${\mathbb S}^3$. 
Let $J$ be a connected component of $\varphi^{-1}(f_n(e_1))$ illustrated in Figure \ref{covering1}.  
Let $l\subset\partial N(k)$ be a meridian of $N(k)$ containing the point $f_n(G_n)\cap \partial N(k)=g_n(G_n)\cap \partial N(k)$. 
Then $l$ is a longitude of $W$. 
In $U$ we find two copies of $f_{n-1}(G_{n-1})$ in $\varphi^{-1}(f_n(G_n)\cup l)$ containing $J$ as illustrated in Figure \ref{covering1}. 
In case $n=1$ we have two copies of a Hopf link. 
Suppose that there is a sequence of crossing changes between $f_n(e_1)$ and $f_n(e_2)$ that deforms $f_n(G_n)$ to $g_n(G_n)$. 
Each of crossing change lifts to crossing changes in $U$. 
In particular a crossing change in $W$ lifts to a crossing change between $J$ and one of the connected components of $\varphi^{-1}(f_n(e_2))$. Since each connected component of $\varphi^{-1}(g_n(e_1))$ bounds a disk in $U\setminus\varphi^{-1}(g_n(G_n\setminus e_1))$ we see that each copy of $f_{n-1}(G_{n-1})$ should be deformed into a copy of $g_{n-1}(G_{n-1})$ by these crossing changes in $U$. 
Therefore the number of crossing changes should be greater than or equal to the twice of $d_{1,2}(f_{n-1},g_{n-1})$. 
See for example Figure \ref{covering2} where the case $n=2$ is illustrated. 
Then we inductively have $d_{1,2}(f_n,g_n)\geq 2^n$. 
$\Box$


\vskip 3mm

\begin{figure}[htbp]
      \begin{center}
\scalebox{0.7}{\includegraphics*{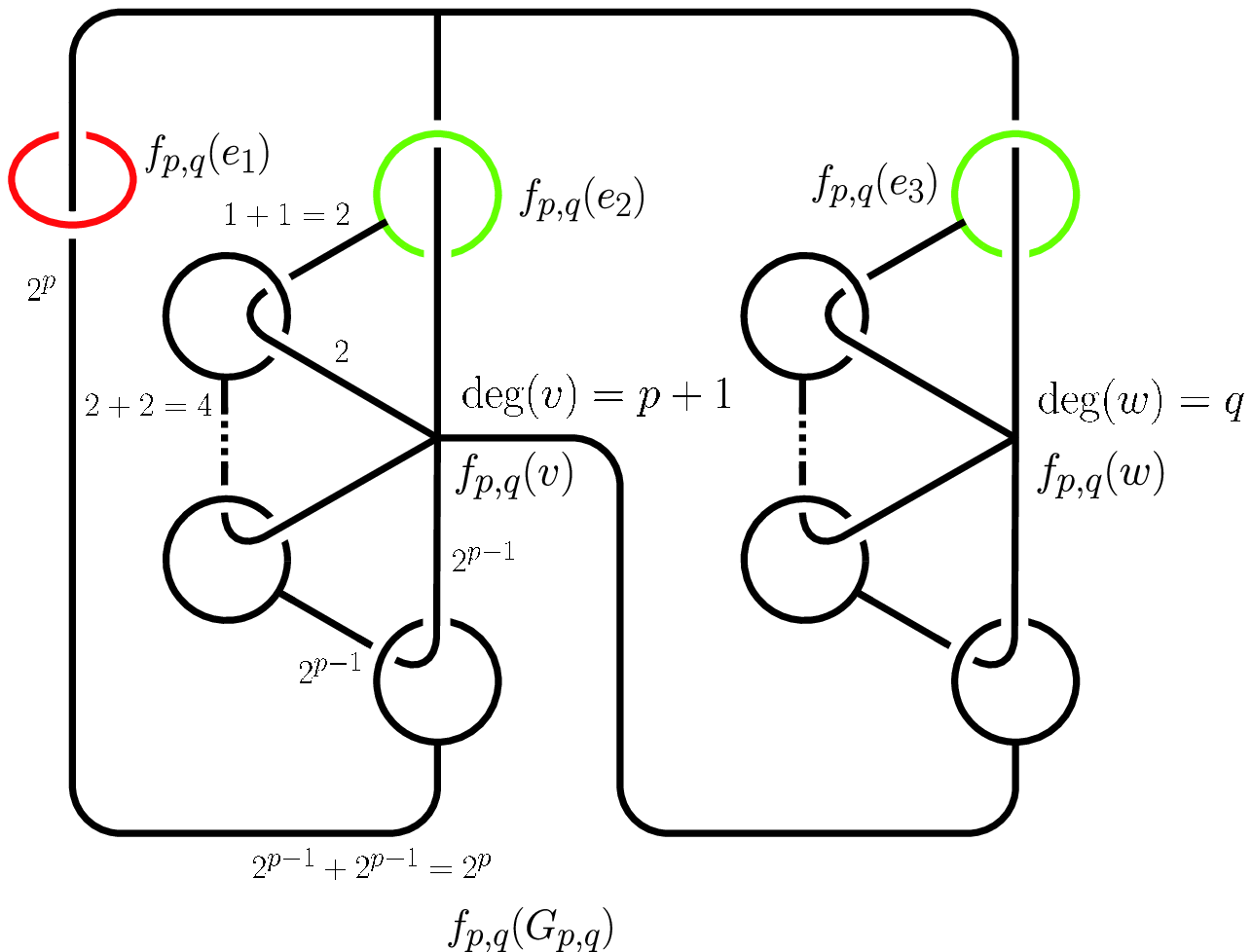}}
      \end{center}
   \caption{}
  \label{pq}
\end{figure}

\begin{Lemma}\label{lemma}
Let $p$ and $q$ be positive integers with $p\leq q$. Let $G_{p,q}$ be a finite graph and $f_{p,q}:G_{p,q}\to{\mathbb S}^3$ a spatial embedding of $G_{p,q}$ as illustrated in Figure \ref{pq}. 
Let $g_{p,q}:G_{p,q}\to{\mathbb S}^3$ be a spatial embedding of $G_{p,q}$ that is different from $f_{p,q}$ only on $e_1$ such that $g_{p,q}(e_1)$ bounds a disk whose interior is disjoint from $g_{p.q}(G_{p,q})$. 
Then $d_{\{[e_1,e_2],[e_1,e_3]\}}(f_{p,q},g_{p,q})=2^p$. 
\end{Lemma}

\vskip 3mm

\begin{figure}[htbp]
      \begin{center}
\scalebox{0.55}{\includegraphics*{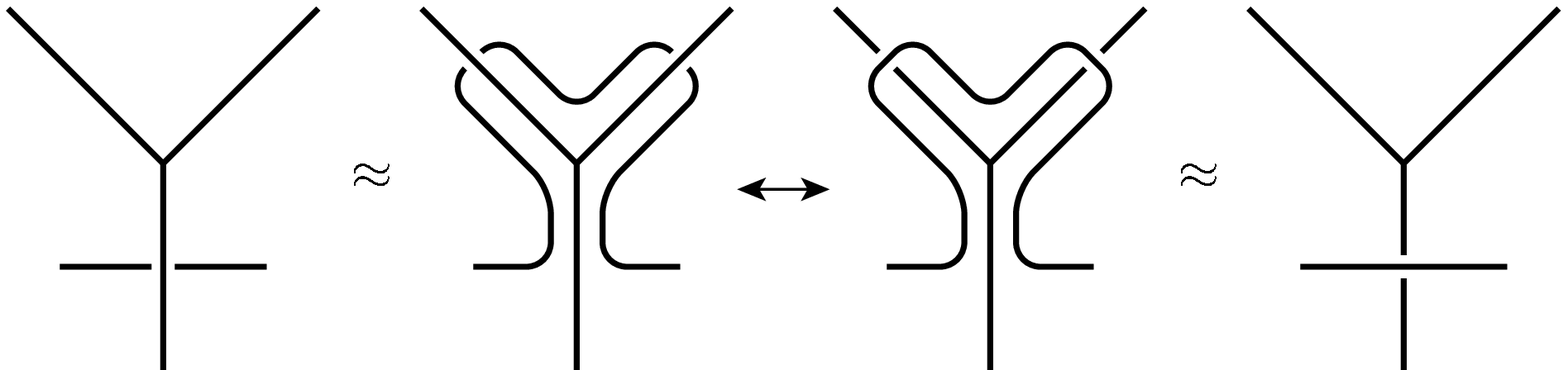}}
      \end{center}
   \caption{}
  \label{replacing2}
\end{figure}

\noindent{\bf Proof.} 
First we show $d_{\{[e_1,e_2],[e_1,e_3]\}}(f_{p,q},g_{p,q})\leq 2^p$. 
We repeatedly use the technic of replacing crossing changes as illustrated in Figure \ref{replacing} and Figure \ref{replacing2}. 
In Figure \ref{pq} the number attached to an edge describes the number of crossing changes with $f_{p,q}(e_2)$ that realize a crossing change with the edge. Then we see $d_{\{[e_1,e_2],[e_1,e_3]\}}(f_{p,q},g_{p,q})\leq 2^p$. 

\begin{figure}[htbp]
      \begin{center}
\scalebox{0.7}{\includegraphics*{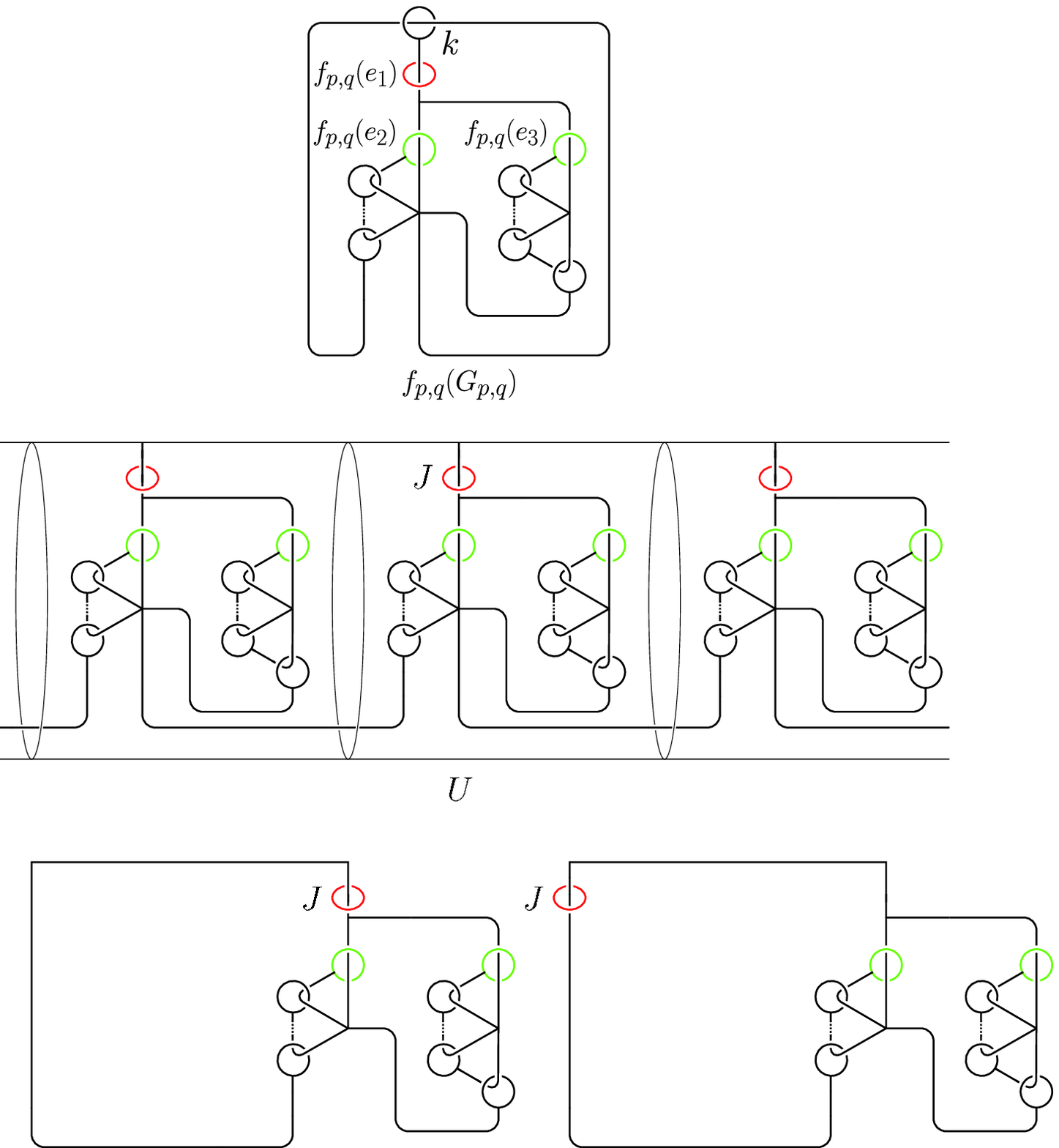}}
      \end{center}
   \caption{}
  \label{pq2}
\end{figure}

Next we show $d_{\{[e_1,e_2],[e_1,e_3]\}}(f_{p,q},g_{p,q})\geq 2^p$ by an induction on $p$. The outline of the proof is similar to that of Theorem \ref{puzzle-1}. 
Suppose that $p=1$. Since $f_{p,q}(G_{p,q})$ is non-splittable \cite{Taniyama2} we have $d_{\{[e_1,e_2],[e_1,e_3]\}}(f_{p,q},g_{p,q})>0$. 
Then by considering the changes of the linking number of $f_{p,q}(e_1)$ and $f_{p,q}(e_2)$ or $f_{p,q}(e_3)$ we have $d_{\{[e_1,e_2],[e_1,e_3]\}}(f_{p,q},g_{p,q})\geq 2$. 
Suppose that $p\geq 2$. Let $k$ be a trivial knot illustrated in Figure \ref{pq2}. 
Let $N(k)$ be a regular neighbourhood of $k$ in ${\mathbb S}^3$. We may suppose that $f_{p,q}(G_{p,q})\cap N(k)=g_{p,q}(G_{p,q})\cap N(k)$ is a simple arc. Since $k$ is trivial $W={\mathbb S}^3\setminus{\rm int}N(k)$ is a solid torus. Let $\varphi:U\to W$ be the universal covering projection. Then we have the preimage $\varphi^{-1}(f_{p,q}(G_{p,q})\cap W)$ in $U$. 
See Figure \ref{pq2}. 
As in the proof of Theorem \ref{puzzle-1} we may suppose that $U$ is a subset of ${\mathbb S}^3$. 
Let $J$ be a connected component of $\varphi^{-1}(f_{p,q}(e_1))$ illustrated in Figure \ref{pq2}. 
Let $l\subset\partial N(k)$ be a meridian of $N(k)$ containing the point $f_{p,q}(G_{p,q})\cap \partial N(k)=g_{p,q}(G_{p,q})\cap \partial N(k)$. 
Then $l$ is a longitude of $W$. 
In $U$ we find two copies of $f_{p-1,q}(G_{p-1,q})$ in $\varphi^{-1}(f_{p,q}(G_{p,q})\cup l)$ containing $J$ as illustrated in Figure \ref{pq2}. 
Therefore as in the proof of Theorem \ref{puzzle-1} we see that $d_{\{[e_1,e_2],[e_1,e_3]\}}(f_{p,q},g_{p,q})$ is greater than or equal to the twice of $d_{\{[e_1,e_2],[e_1,e_3]\}}(f_{p-1,q},g_{p-1,q})$. 
Then we inductively have $d_{\{[e_1,e_2],[e_1,e_3]\}}(f_{p,q},g_{p,q})\geq 2^p$. 
$\Box$


\vskip 3mm

\begin{figure}[htbp]
      \begin{center}
\scalebox{0.8}{\includegraphics*{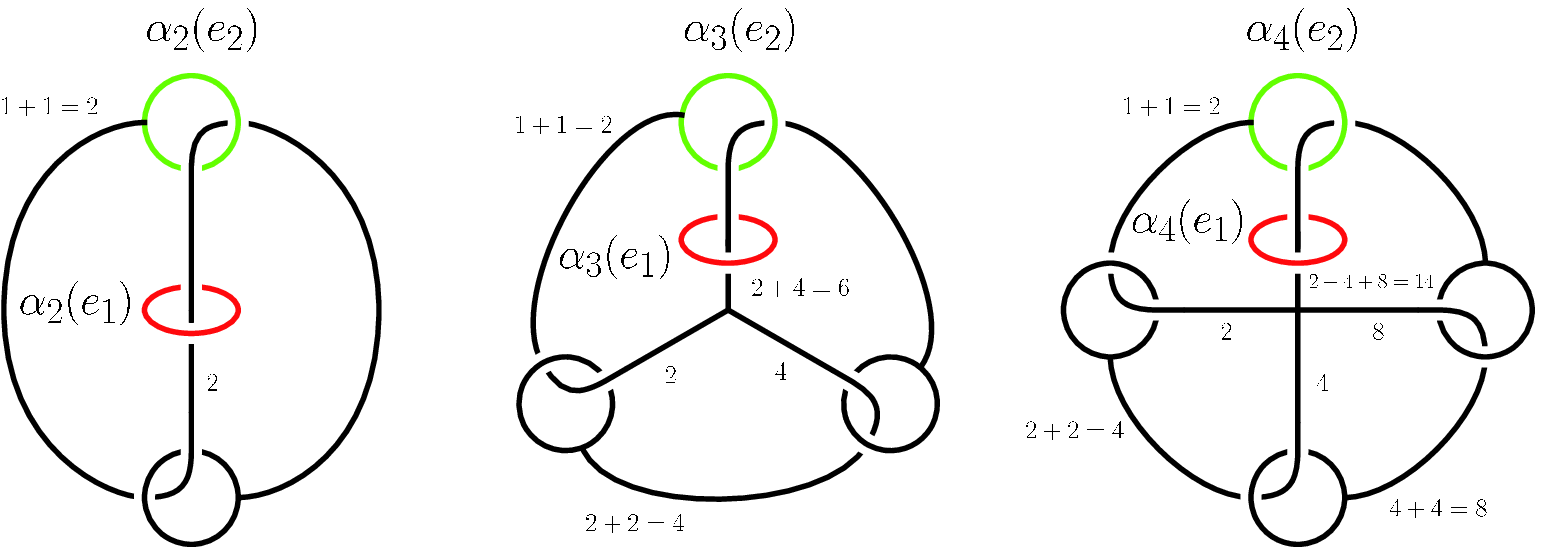}}
      \end{center}
   \caption{}
  \label{puzzle2-2}
\end{figure}

\noindent{\bf Proof of Theorem \ref{puzzle-2}.} 
First we show $d_{1,2}(\alpha_n,\beta_n)\leq2^n-2$. 
We repeatedly use the technic illustrated in Figure \ref{replacing} and Figure \ref{replacing2}. 
In Figure \ref{puzzle2-2} the number attached to an edge describes the number of crossing changes with $\alpha_n(e_2)$ that realize a crossing change with the edge. Then we see $d_{1,2}(\alpha_n,\beta_n)\leq 2+4+8+\cdots+2^{n-1}=2^n-2$.

Next we show $d_{1,2}(\alpha_n,\beta_n)\geq2^n-2$. The outline of the proof is similar to that of Theorem \ref{puzzle-1}. 

\begin{figure}[htbp]
      \begin{center}
\scalebox{0.6}{\includegraphics*{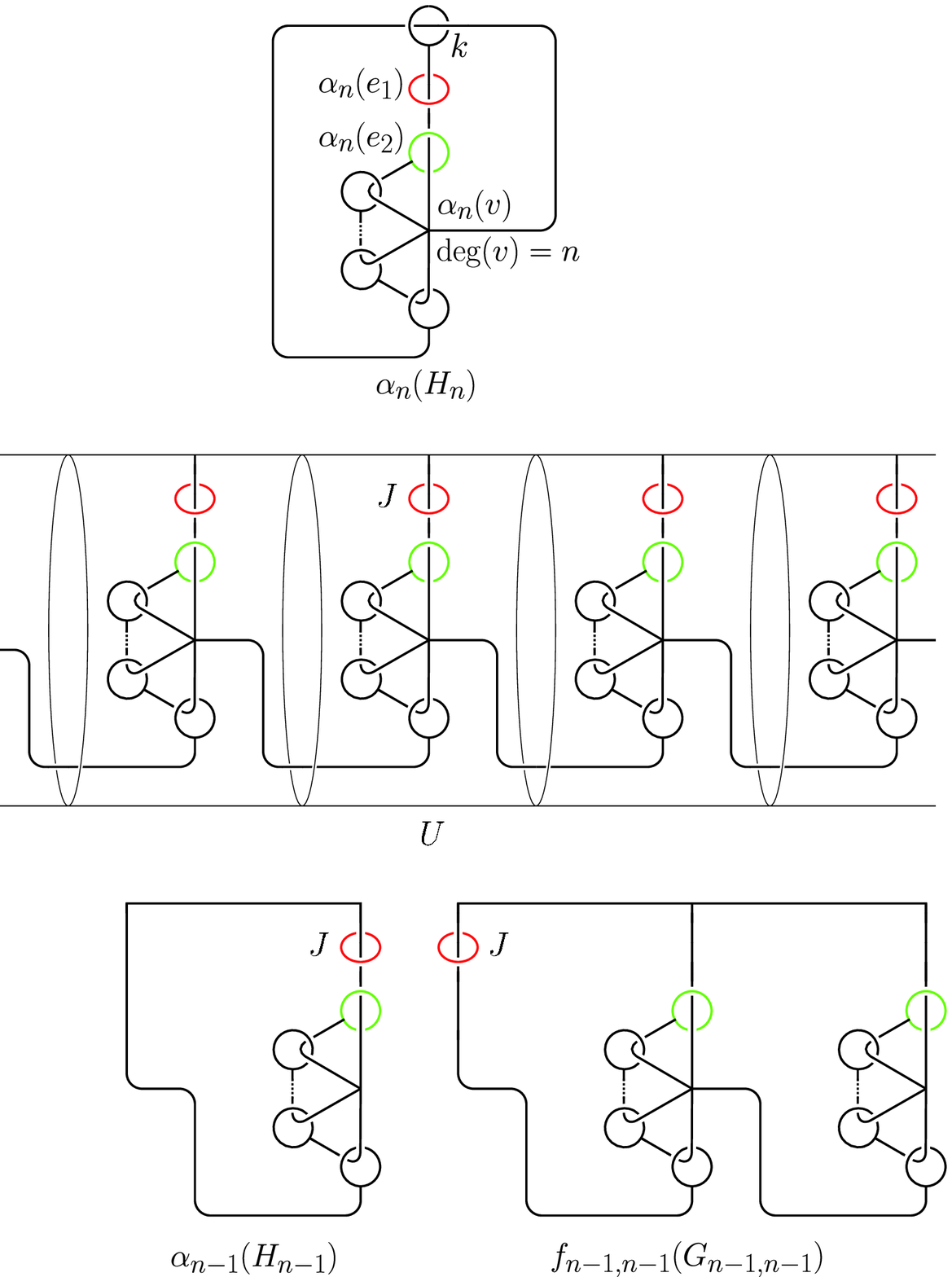}}
      \end{center}
   \caption{}
  \label{puzzle4}
\end{figure}

Suppose that $n=2$. Since $\alpha_2(H_2)$ is non-splittable \cite{Taniyama2} we have $d_{1,2}(\alpha_2,\beta_2)>0$. 
Then by considering the changes of the linking number of $\alpha_2(e_1)$ and $\alpha_2(e_2)$ we have $d_{1,2}(\alpha_2,\beta_2)\geq 2$. 
Suppose that $n\geq 3$. 
Let $k$ be a trivial knot illustrated in Figure \ref{puzzle4} and $N(k)$ a regular neighbourhood of $k$ in ${\mathbb S}^3$ such that $\alpha_{n}(H_{n})\cap N(k)=\beta_{n}(H_{n})\cap N(k)$ is a simple arc. Let $W={\mathbb S}^3\setminus{\rm int}N(k)$ be a solid torus. Let $\varphi:U\to W$ be the universal covering projection. Then we have the preimage $\varphi^{-1}(\alpha_{n}(H_{n})\cap W)$ in $U$. 
See Figure \ref{puzzle4}. 
As in the proof of Theorem \ref{puzzle-1} we may suppose that $U$ is a subset of ${\mathbb S}^3$. 
Let $J$ be a connected component of $\varphi^{-1}(\alpha_{n}(e_1))$ illustrated in Figure \ref{puzzle4}. 
Let $l\subset\partial N(k)$ be a meridian of $N(k)$ containing the point $\alpha_{n}(H_{n})\cap \partial N(k)=\beta_{n}(H_{n})\cap \partial N(k)$. 
Then $l$ is a longitude of $W$. 
In $\varphi^{-1}(\alpha_{n}(H_{n})\cup l)$ in $U$ we find a copy of $\alpha_{n-1}(H_{n-1})$ and a copy of $f_{n-1,n-1}(G_{n-1,n-1})$ both containing $J$ as illustrated in Figure \ref{puzzle4}. 
Therefore as in the proof of Theorem \ref{puzzle-1} we see that $d_{1,2}(\alpha_n,\beta_n)$ is greater than or equal to the sum of $d_{1,2}(\alpha_{n-1},\beta_{n-1})$ and $d_{\{[e_1,e_2],[e_1,e_3]\}}(f_{n-1,n-1},g_{n-1,n-1})$. 
Since $(2^{n-1}-2)+2^{n-1}=2^n-2$ we inductively have $d_{1,2}(\alpha_n,\beta_n)\geq2^n-2$. 
$\Box$

\vskip 3mm

\begin{Proposition}\label{well-definedness}
The link type of $l_n$ is independent of the choice of the defining edge $e$ of $T$. 
\end{Proposition}

\vskip 3mm

\begin{figure}[htbp]
      \begin{center}
\scalebox{0.5}{\includegraphics*{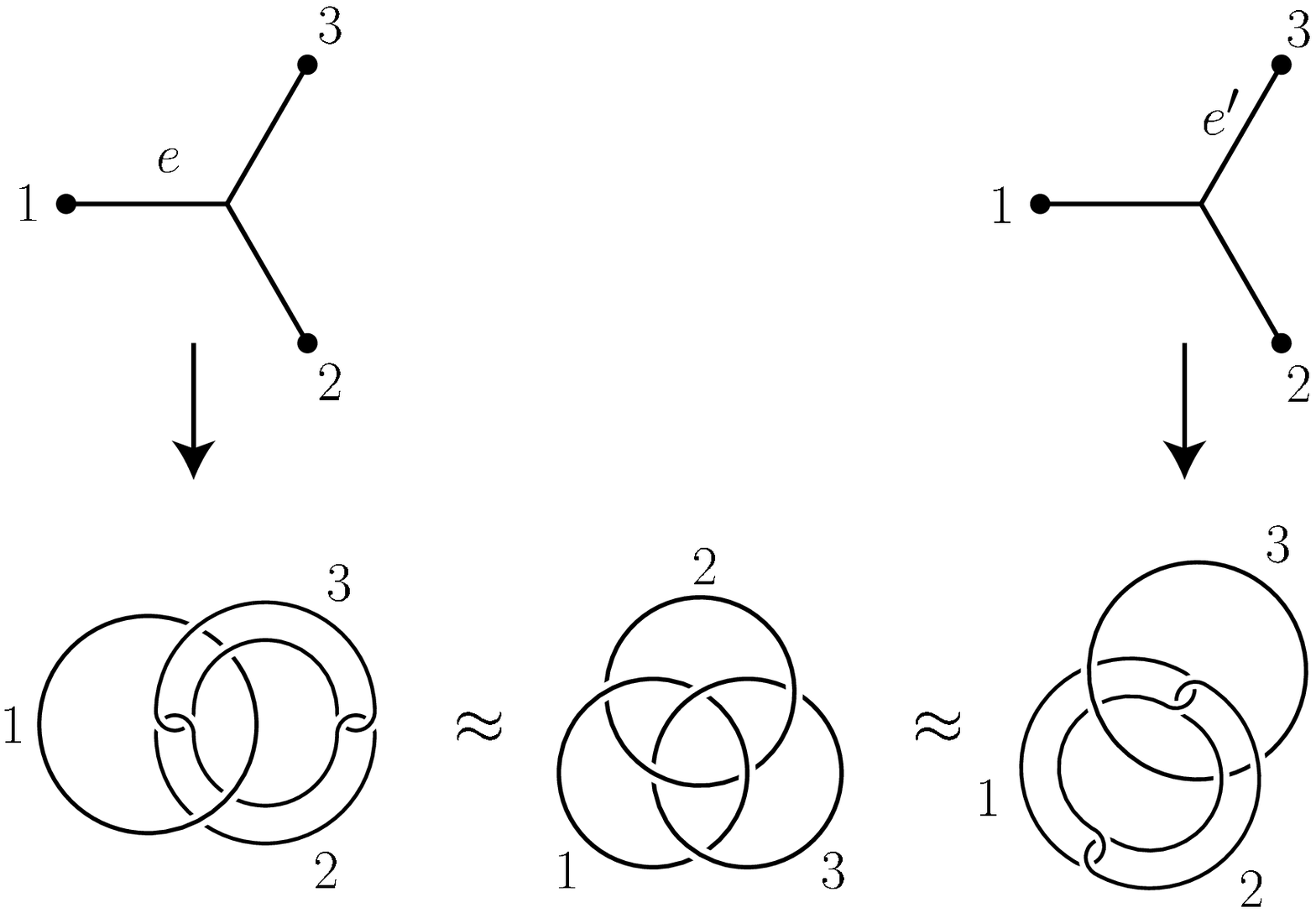}}
      \end{center}
   \caption{}
  \label{Borromean}
\end{figure}

\noindent{\bf Proof.} It is sufficient to show that the link type of $l_n$ is unchanged when we replace $e$ by an edge $e'$ adjacent to $e$. This follows from the fact that the Borromean rings is a Bing double of a Hopf link and the Borromean rings has $3$-symmetry with respect to the components. See Figure \ref{Borromean}. $\Box$

\vskip 3mm

\begin{Lemma}\label{tree}
Let $T$ be a uni-trivalent tree and $v,w,x$ mutually distinct vertices in $V_1(T)$. Let $P(v,w)$ be the path in $T$ with $\partial P(v,w)=\{v,w\}$. 
Let $y$ be the vertex of $T$ adjacent to $x$. Let $s$ and $t$ be the other vertices of $T$ adjacent to $y$. 
Let $T'$ be a uni-trivalent tree obtained from $T$ by deleting $x$, $y$ and the edges incident to $y$ and by joining $s$ and $t$ by a new edge $e'$. 
Let $a$ and $b$ {\rm (}resp. $a'$ and $b'${\rm )} be the components of $L_T$ {\rm (}resp. $L_{T'}${\rm )} corresponding to $v$ and $w$ respectively. 
Let ${\mathcal F}=\{[a,b]\}$ and ${\mathcal F'}=\{[a',b']\}$. 
If $y$ is not a vertex of $P(v,w)$, then $u_{\mathcal F}(L_T)\geq u_{\mathcal F'}(L_{T'})$. 
If $y$ is a vertex of $P(v,w)$, then $u_{\mathcal F}(L_T)\geq 2u_{\mathcal F'}(L_{T'})$. 
\end{Lemma}

\vskip 3mm

\begin{figure}[htbp]
      \begin{center}
\scalebox{0.5}{\includegraphics*{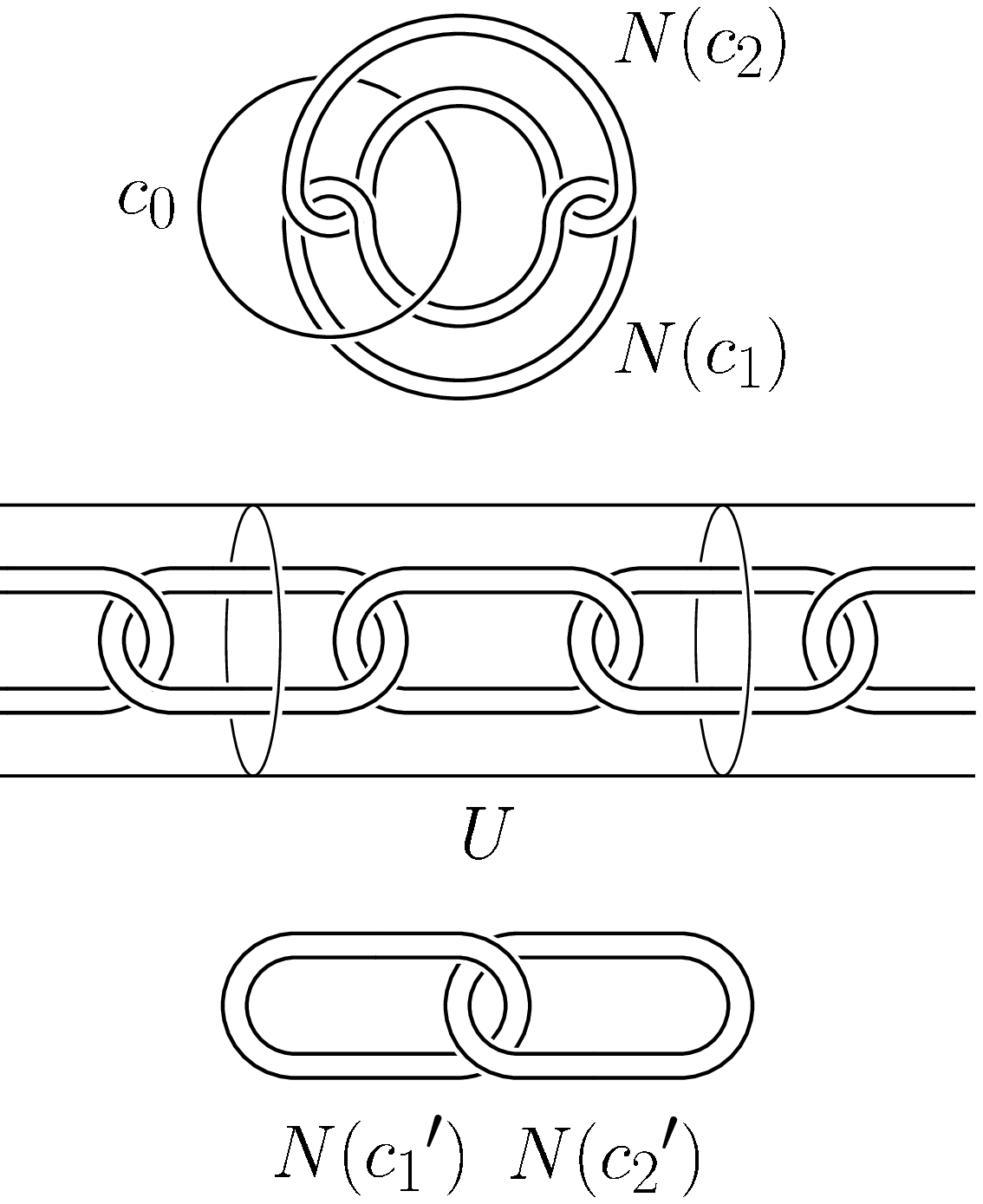}}
      \end{center}
   \caption{}
  \label{Borromean2}
\end{figure}

\noindent{\bf Proof.} 
Let $c_0$ be a component of $L_T$ corresponding to $x$. Let $e$ be the edge of $T$ joining $x$ and $y$. By Proposition \ref{well-definedness} we may suppose that $e$ is the defining edge of $L_T$. Let $M=c_0\cup c_1\cup c_2$ be Borromean Rings and $N(c_0)$, $N(c_1)$ and $N(c_2)$ mutually disjoint regular neighbourhoods of $c_0$, $c_1$ and $c_2$ respectively such that $L_T$ is contained in $c_0\cup N(c_1)\cup N(c_2)$ and obtained from $M$ by repeatedly taking Bing doubles. 
If $y$ is not a vertex of $P(v,w)$ then $a\cup b$ is contained in $N(c_1)$ or $N(c_2)$. 
If $y$ is a vertex of $P(v,w)$ then $a\cup b$ is not contained in $N(c_1)$ nor $N(c_2)$. 
We note that $W={\mathbb S}^3\setminus{\rm int}N(c_0)$ is a solid torus containing $N(c_1)\cup N(c_2)$. 
Let $\varphi:U\to W$ be the universal covering projection. 
Then we have the preimage $\varphi^{-1}(N(c_1)\cup N(c_2))$ in $U$ and $\varphi^{-1}(L_T\setminus c_0)$ in $\varphi^{-1}(N(c_1)\cup N(c_2))$. 
By Proposition \ref{well-definedness} we may suppose that $e'$ is the defining edge of $L_{T'}$. 
Then there is a Hopf link $H={c_1}'\cup {c_2}'$ and mutually disjoint regular neighbourhoods $N({c_1}')$ and $N({c_2}')$ of ${c_1}'$ and ${c_2}'$ respectively such that $L_{T'}$ is contained in $N({c_1}')\cup N({c_2}')$ and obtained from $H$ by repeatedly taking Bing doubles. We may also suppose that the pair $(N(c_1),L_T\cap N(c_1))$ is homeomorphic to the pair $(N({c_1}'),L_{T'}\cap N({c_1}'))$ and the pair $(N(c_2),L_T\cap N(c_2))$ is homeomorphic to the pair $(N({c_2}'),L_{T'}\cap N({c_2}'))$. 
Therefore there exist countably many copies of $L_{T'}$ in $\varphi^{-1}(L_T\setminus c_0)$. See Figure \ref{Borromean2}. 
Then we have $u_{\mathcal F}(L_T)\geq u_{\mathcal F'}(L_{T'})$. 
If $y$ is a vertex of $P(u,v)$, then $a\cup b$ is not contained in $N(c_1)$ nor $N(c_2)$ and the situation is lifted to $U$. 
Then by the same reason as in the proofs of Theorem \ref{puzzle-1} and Theorem \ref{puzzle-2} we have $u_{\mathcal F}(L_T)\geq 2u_{\mathcal F'}(L_{T'})$. 
$\Box$

\vskip 3mm

\noindent{\bf Proof of Theorem \ref{Milnor}.} 
Let $v$ and $w$ be the vertices in $V_1(T)$ corresponding to $a$ and $b$ respectively. Let $P(v,w)$ be the path in $T$ with $\partial P(v,w)=\{v,w\}$. 

\begin{figure}[htbp]
      \begin{center}
\scalebox{0.5}{\includegraphics*{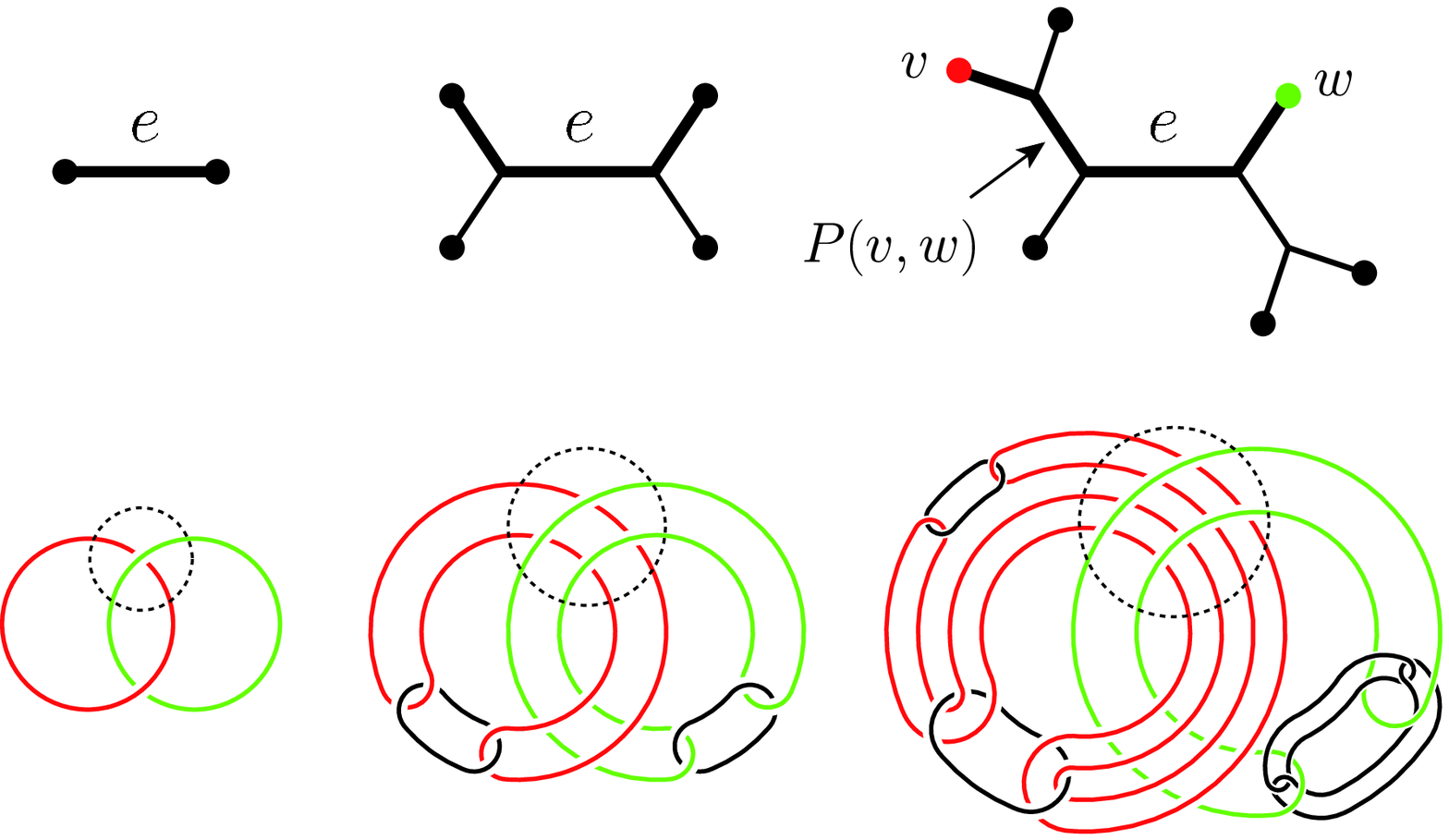}}
      \end{center}
   \caption{}
  \label{Milnor2}
\end{figure}

First we show $u_{\mathcal F}(L_T)\leq 2^{d(a,b)-1}$. 
Let $e$ be an edge of $P(v,w)$. 
By Proposition \ref{well-definedness} we may suppose that $e$ is the defining edge of $L_T$. 
In the process of taking a Bing double involving a component corresponding a vertex of $P(v,w)$ we may suppose that the other component is small enough so that it has only four crossings in a diagram. Thus we have a diagram of $L_T$ and $2^{d(a,b)-1}$ crossings of it such that changing all of them will turn $L_T$ into a trivial link. See for example Figure \ref{Milnor2}. 

Next we show $u_{\mathcal F}(L_T)\geq 2^{d(a,b)-1}$ by an induction on the number of the edges of $T$. 
If $T$ has exactly one edge, then $L_T$ is a Hopf link and $d(a,b)=1$. Then we have $u_{\mathcal F}(L_T)=1=2^{d(a,b)-1}$. 
Suppose that $T$ has three or more edges. Assume that the claim is true for all trees of fewer edges. 
Let $x$ be a vertex in $V_1(T)\setminus\{v,w\}$. Let $y$ be the vertex of $T$ adjacent to $x$. 
Let $T'$ be the graph obtained from $T$ as described in Lemma \ref{tree}. 
Let $a'$ and $b'$ be the components of $L_{T'}$ corresponding to $v$ and $w$ respectively. 
Let ${\mathcal F'}=\{[a',b']\}$. 
Let $P'(v,w)$ be the path in $T'$ with $\partial P'(v,w)=\{v,w\}$. 
If $y$ is not a vertex of $P(v,w)$, then $P(v,w)=P'(v,w)$ and the number of the edges of $P'(v,w)$ is equal to that of $P(v,w)$. 
Therefore $d(a,b)=d(a',b')$. 
By Lemma \ref{tree} we have $u_{\mathcal F}(L_T)\geq u_{\mathcal F'}(L_{T'})$. 
By the assumption $u_{\mathcal F'}(L_{T'})\geq 2^{d(a',b')-1}$. 
Thus we have $u_{\mathcal F}(L_T)\geq 2^{d(a,b)-1}$. 
If $y$ is a vertex of $P(v,w)$, then the number of the edges of $P'(v,w)$ is one less than that of $P(v,w)$. 
Therefore $d(a,b)=d(a',b')+1$. 
By Lemma \ref{tree} we have $u_{\mathcal F}(L_T)\geq 2u_{\mathcal F'}(L_{T'})$. 
By the assumption $u_{\mathcal F'}(L_{T'})\geq 2^{d(a',b')-1}$. 
Thus we have $u_{\mathcal F}(L_T)\geq 2^{d(a,b)-1}$. 
$\Box$

\vskip 3mm

\section*{Acknowledgments} In 2012 a Japanese company East Entertainment Inc. asked the author to provide a problem in knot theory for a TV program of Fuji Television Network, Inc. entitled Takeshi-no-Komadai-Sugakuka. It was an entertainment program based on mathematics. Then the author considered puzzle ring problem and he have found a new proof of Theorem \ref{puzzle-1}. The author is grateful to East Entertainment Inc. for giving him such an opportunity.

\vskip 3mm

{\normalsize
\end{document}